\documentclass[10]{article}

\def\ba{\mathbf{a}}   
\def\bb{\mathbf{b}}   
\def\be{\mathbf{e}}   
\def\bm{\mathbf{m}}   
\def\bx{\mathbf{x}}   
\def\by{\mathbf{y}}   

\def\bA{\mathbf{A}}   
\def\bB{\mathbf{B}}   
\def\bJ{\mathbf{J}}   

\def\bX{\mathbf{X}}   
\def\bY{\mathbf{Y}}   
\def\G{\mathbb{G}} 

\def\R{\mathbb{R}}  

\def\no{\noindent}

\def\ul{\underline}
\def\ol{\overline}

\def\beq{\begin{equation}}
\def\eeq{\end{equation}}
\def\con{\overline}

\def\w{\wedge}

\usepackage{array}
\usepackage{tabularx}
\usepackage{amsfonts}
\usepackage{amssymb}
\usepackage{tikz-cd}
\usetikzlibrary{matrix,arrows,decorations.pathmorphing}

\def\bpm{\begin{pmatrix}}
	\def\epm{\end{pmatrix}}
\makeindex

\title{Spheroidal Domains and Geometric Analysis in Euclidean Space} 
\author{Garret Sobczyk \\  
Universidad de las Am\'ericas-Puebla \\
 Departamento de Actuaría F\'isica y Matem\'aticas \\
72820 Puebla, Pue., M\'exico}

\begin{document}
\maketitle
\begin{abstract} 
Clifford's {\it geometric algebra} has enjoyed phenomenal development over the last 60 years by mathematicians, theoretical physicists, engineers and computer scientists in robotics, artificial intelligence and data analysis, introducing a myriad of different and often confusing notations. The geometric algebra of Euclidean 3-space, the natural generalization of both the well-known Gibbs-Heaviside vector algebra, and Hamilton's quaternions, is used here to study spheroidal domains, spheroidal-graphic projections, the Laplace equation and its Lie algebra of symmetries. The Cauchy-Kovalevska extension and the Cauchy kernel function are treated in a unified way. The concept of a {\it quasi-monogenic} family of functions is introduced and studied.

 \smallskip
 
\no {\em AMS Subject Classification:} 15A66, 30A05, 35J05.

\smallskip

\no {\em Keywords: geometric analysis, Clifford analysis, spheroidal Laplacian, quasi-monogenic functions.}
\end{abstract}

 \section*{0\quad Introduction}
 
 Two main scientific communities utilizing William Kingdon Clifford's {\it geometric algebra} have been in development over the last 60 or more years. The {\it Clifford analysis} community has developed Clifford algebra primarily as the natural generalization to higher dimensions of the ubiquitous complex analysis of analytic functions, which underlies much of modern mathematics and theoretical physics. The second community, which I dub the {\it geometric analysis} community, has stressed the more general development of geometric algebra as the natural generalization of the real number system to include the concept of direction. The first community consists in large part of mathematicians, where as the second community consists of a more diverse group of people in mathematics, theoretical physics and computer scientists, and engineers interested in diverse applications such as robotics, artificial intelligence and data analysis. Geometric algebra $\G_3$ is the natural extension of the popular Gibbs-Heaviside vector algebra still universally employed by many engineers and scientists today. 
 
 Whereas there is a great deal of overlap between these groups, namely the usage of Clifford algebra, invented by W.K. Clifford in the years shortly before his death in 1879, the different symbolisms and notations employed has lead to a general lack of communications between the two groups. It is the belief of the present author that a greater communication between the two groups would be advantageous to both groups. {\it Spheroidal domains}, usually studied in terms of {\it quaternion analysis}, are here reformulated in the {\it geometric analysis} of Euclidean space. Spherical domains and spherical harmonics are a limiting case of spheroidal domains and spheroidal harmonics \cite{BBS1997}. 
 
 Section 1, sets down the basic definitions of prolate and oblate spheroidal coordinates in terms of the associative geometric algebra $\G_3$ of Euclidean space $\R^3$,     
 \[ \G_3:=\G(\R^3 )=\R[\be_0,\be_1,\be_2  ],   \]
where $\be_k$ are three orthogonal {\it anti-commuting} unit vectors along the respective $x_k$-axis for $k=0,1,2$. That is
\[  \be_k^2 =1, \quad {\rm and} \quad \be_{jk}:=\be_j \be_k =-\be_k \be_j=- \be_{kj},        \]
for $k\ne j$. The notation used is meant to suggest that the real number system $\R$ is extended to include the three unit orthogonal vectors $\be_k$ and their geometric sums and products \cite{H/S, Sob2019, SNF}. As seen in later sections, spheroidal coordinates find their importance in being one of 11 orthogonal separable coordinate systems, \cite[p.\,40]{BKM1976}.  
 
Section 2, studies {\it spheroidal-graphic projection} of the unit prolate and oblate spheroids onto the two dimensional plane, the natural generalization of more famous {\it stereographic projection}. This serves to help unfamiliar readers come to grips with the concept of prolate and oblate spheroids, which may be otherwise unfamiliar to them.

Sections 3 and 4, introduce prolate and oblate spheroidal gradients and Laplacians in a unified way, taking advantage of the rich geometric structure of the geometric algebra $\G_3$.  

Section 5, studies solutions of the Laplace equation, in both the prolate and oblate cases, using the well-known method of {\it separation of variables}.

Section 6, briefly considers the beautiful theory of the Lie algebra of symmetry operators, which gives insight into the century long history of the subject. 

Section 7, shows how Clifford analysis can be incorporated directly into the body of the more comprehensive geometric analysis, unifying the otherwise different approaches. The concept of a quaternion arises naturally in the even sub-algebra of the geometric algebra $\G_3$ of Euclidean 3-space.  
 As an application, the Cauchy kernel function is used to generate a monogenic hypercomplex power series, \cite[(3.6),(3.7)]{CFM2017}. The {\it Cauchy-Kovalevska} extension, a method for generating higher order monogenic functions, has been treated by many authors, \cite{CFM2017} and \cite[p.151]{DSS1992}. By using a simple idea suggested by this extension, a family of curl-free {\it quasi-mononogenic} functions is generated.

In an Appendix, a Mathematica Package is included giving solutions to the separable differential equations explored in Section 5.   

    \section{Prolate and oblate spheroidal coordinates}

    Let $\G_3:=\R(e_0,e_1,e_2 )$ be the geometric algebra of $3$-dimensional
   Euclidean space $\R^3$. The position vectors $\bx$ and $\by$ in {\it prolate and oblate spheroidal coordinates} $(\eta,\theta, \varphi) $ can be defined, respectively, in terms of the {\it complex-number} like
   quantity 
   \beq z:=\frac{1}{2}\Big( e^{\eta + \theta I_p}+ e^{-(\eta +\theta  I_p)}\Big)=\cosh(\eta+I_p \theta)=\cosh \eta \cos \theta +I_p \sinh \eta \sin \theta \label{complexroid} \eeq
  where $I_p:=e_p e_0$ has square minus one for $e_p = \be_1\cos \varphi  +  \be_2\sin\varphi$, and where $ \eta \ge 0 ,\ \mu>0 , \ \varphi \in [0,2\pi)$, 
 $ \theta \in [0,\pi]$.
 
  For $\bx$, 
 \beq  \bx:=x_0 e_0+x_p e_p =\mu ze_0 = \mu e_0 \ol z = \mu e_0 \cosh(\eta-I_p \theta) \label{spheroidalunited} \eeq
  where
   \[  x_0 =\mu \cosh\eta \cos \theta ,\ \ x_p := \sqrt{x_1^2+x_2^2}=\mu\sinh \eta \sin \theta, \]
     \[  x_1 =\mu \cos \varphi \sinh\eta \sin \theta  ,\ \ x_2 := \mu \sin \varphi\sinh \eta \sin \theta, \]
  in the prolate case, and
  \beq \by:=y_0 e_0+y_p e_p =\mu z_\eta e_0 =\mu e_0\ol z_\eta= \mu \sinh(\eta+I_p \theta)e_0 \label{spheroidalunitedo} \eeq
  where $z_\eta := \partial_{\eta}z$, so that
   \[  y_0 =\mu\sinh \eta \cos \theta  ,\ \ y_p := \sqrt{y_1^2+y_2^2}=\mu\cosh\eta \sin \theta ,\]
    \[  y_1 =\mu \cos \varphi \cosh\eta \sin \theta  ,\ \ y_2 := \mu \sin \varphi\cosh \eta \sin \theta, \]
   in the oblate case, \cite{BKM1976,Hob1931,wofram}.\footnote{Different conventions are used for oblate coordinates. The oblate coordinates used here are the same as in \cite{BKM1976,Hob1931}. }   
   
   Equations (\ref{spheroidalunited}) and  (\ref{spheroidalunitedo}) give a direct relationship between prolate and oblate coordinates, and their expression in terms of the {\it quaternion-like} quantities $z$ and $\ol z$. 
  Since the bivector $I_p=e_p e_0$ has square $-1$, it behaves the same as the imaginary unit $i=\sqrt{-1}$. Note that 
   \[ \dot I_p:=\partial_\varphi I_p =\dot e_p e_0=(- \be_1\sin \varphi  +  \be_2\cos\varphi)e_0 \]
   also has square $-1$, as does the quantity $I_p\dot I_p=\dot e_p e_p$. Indeed, the bivectors $I_p,J_p:=\dot I_p,K_p=\dot e_pe_p$ obey exactly the same rules as Hamilton's quaternions. The dot over a variable is used to denote the partial derivative with respect to $\varphi$. Thus $\dot z :=z_\varphi=\partial_\varphi z$.  
   
    We also calculate
   \[   \bx^2 = \mu^2 z \ol z =\frac{1}{2}(\cosh 2\eta+\cos 2 \theta), \ \  \by^2=z_\eta \ol z_\eta =\frac{1}{2}(\cosh 2\eta-\cos 2 \theta)  , \] 
   and define the quantities
   \beq \omega_x := |\bx+\mu e_0|+|\bx-\mu e_0 | =\sqrt{(x_0+\mu)^2+x_p^2}+\sqrt{(x_0-\mu)^2+x_p^2} =2\mu \cosh \eta \label{pomega} \eeq
 and
  \beq \ol \omega_x := |\bx + \mu e_0| -|\bx - \mu e_0 | =\sqrt{(x_0+\mu)^2+x_p^2}-\sqrt{(x_0-\mu)^2+x_p^2} =2 \mu \cos \theta, \label{pcomega} \eeq
  in the prolate case.  In the oblate case, 
    \beq \omega_y := |\by + \mu e_p| +|\by - \mu e_p |=\sqrt{\mu^2+y^2 + 2 \mu y_p}+\sqrt{\mu^2+y^2 - 2 \mu y_p}=2 \mu \cosh \eta  \label{obomega} \eeq
  and
  \beq \ol \omega_y := |\by + \mu e_p| -|\by - \mu e_p |=\sqrt{\mu^2+y^2 + 2 \mu y_p}-\sqrt{\mu^2+y^2 - 2 \mu y_p} =2 \mu \sin \theta.  \label{obcomega} \eeq
  
\begin{figure}[h]
 	\begin{center}  	\includegraphics[width=0.50\linewidth]{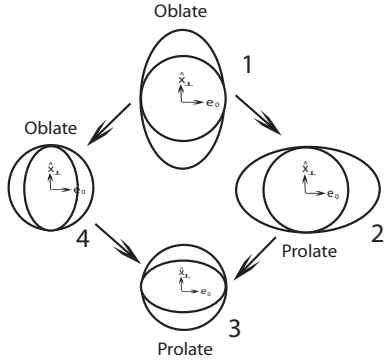}  \end{center}
 	\caption{Of the four unit bounding spheroids pictured, two are oblate and two are prolate, and are rotated around the $\be_0$-axis.  For $\nu \ge 0$, $e^{2 \nu}-\mu^2 =1$ for Cases 1, 2, and
 		$e^{-\nu}+\mu^2 = 1$ for Cases 3, 4, respectively.} 
 	\label{prolateoblate}
 \end{figure}

   The proofs of the equations (\ref{pomega}) - (\ref{obcomega}) are very similar. For the prolate case,
  \[ |\bx \pm \mu e_0|^2 = \mu^2|z \pm 1|^2 =(\cosh\eta \pm \cos \theta )^2 ,    \]
  and for the oblate case,
    \[ |\by \pm \mu e_p|^2 = \mu^2|z_\eta \pm I_p|^2 =(\cosh\eta \pm \sin \theta )^2 .    \]
  Geometrically, $\omega_x$ defined in (\ref{pomega}) and $\omega_y$ defined in (\ref{obomega}) are distances on the bounding unit prolate and oblate spheroids between the focal points located at the points $(0,\pm\mu,0)$ in the prolate cases 2 \& 3, and the focal points located at the points $(0,0,\pm \mu)$ in the oblate cases 4 \& 1in Figure \ref{prolateoblate}, respectively. Similarly, $\con \omega_x$ and $\con \omega_y$ are the distances between the foci of the bounding unit hyperbolic spheroids in the prolate and oblate cases, respectively,
  \cite{wofram}. 
  
  Since $\omega_x=\omega_y$ and $\ol \omega_x(\theta) = \ol \omega_y(\theta + \frac{\pi}{2})$, it follows that
  \[  \omega_x= \sqrt 2 \Big((x^2+\mu^2)+\sqrt{(x^2+ \mu^2)^2-4 \mu^2x_0^2}\Big)^{\frac{1}{2}} \]
   \beq = \sqrt 2 \Big((y^2+\mu^2)+\sqrt{(y^2+ \mu^2)^2-4 \mu^2y_p^2}\Big)^{\frac{1}{2}}=\omega_y ,  \label{omegaxy} \eeq
   and 
   \[  \ol\omega_x(\theta)= \sqrt 2 \Big((x^2+\mu^2)-\sqrt{(x^2+ \mu^2)^2-4 \mu^2x_0^2}\Big)^{\frac{1}{2}} \]
  \beq = \sqrt 2 \Big((y^2+\mu^2)-\sqrt{(y^2+ \mu^2)^2-4 \mu^2y_p^2}\Big)^{\frac{1}{2}}=\ol\omega_y(\theta + \frac{\pi}{2}) .  \label{omegabxy} \eeq
  
  The equations (\ref{omegaxy}) and (\ref{omegabxy}) define a set of four bounding unit spheroids, pictured in Figure \ref{prolateoblate}.
  \begin{itemize}
  	\item[3.] $\frac{\cosh[\eta + I \theta] }{\cosh \eta}e_0 =\Big(\cos \theta +I_p \tanh \eta \sin \theta\Big)e_0, \, \mu \cosh \eta =1  \iff  e^{-\nu} := \tanh \eta $
  	\item[2.] $\frac{\cosh[\eta + I \theta] }{\sinh \eta}e_0 =\Big(\coth \eta \cos \theta +I_p  \sin \theta\Big)e_0, \, \mu \cosh \eta =1  \iff  e^{\nu} := \coth \eta $ 
  	\item[4.] $\frac{\sinh[\eta + I \theta] }{\cosh \eta}e_0 =\Big(\tanh \eta \cos \theta +I_p  \sin \theta\Big)e_0, \, \mu \cosh \eta =1  \iff  e^{-\nu} := \tanh \eta $ 
  	\item[1.] $\frac{\sinh[\eta + I \theta] }{\sinh \eta}e_0 =\Big(\cos \theta +I_p \coth \eta \sin \theta\Big)e_0, \, \mu \cosh \eta =1  \iff  e^{\nu} := \coth \eta $ 	
  \end{itemize}
  
For the coordinates $(\eta, \theta, \varphi)$, the partial derivatives
\[ z_\eta :=\partial_\eta z = \sinh(\eta + I_p \theta), \ z_\theta :=\partial_\theta z = I_p\sinh(\eta + I_p \theta),   \]
\[ z_{\eta\eta} :=\partial_\eta^2 z =z, \ z_{\theta \theta} :=\partial_\theta^2 z =-z , \ z_{\eta \theta} := \partial_{\eta}\partial_{\theta}z=I_p z,   \]
and 
\[ z_\varphi  = \partial_\varphi z=\dot z =\dot I_p \sinh \eta \sin \theta = \dot I_p(e_p\cdot \bx) , \  z_{\varphi \varphi}  = \partial_\varphi^2 z=-I_p \sinh \eta \sin \theta, \]
 and are use to calculate,
\[ \bx_\eta := \partial_\eta \bx = \mu z_\eta e_0, \   \bx_\theta := \partial_\theta \bx = \mu z_\theta e_0, \  \bx_\varphi := \partial_\varphi \bx = \mu z_\varphi e_0=\mu \dot e_p \sinh \eta \sin \theta \] 
for the {\it prolate orthogonal tangent vectors} $\{\bx_\eta,\bx_\theta, \bx_\varphi \}$. The corresponding {\it orthogonal reciprocal frame} $\{\bx^\eta,\bx^\theta, \bx^\varphi \}$ is defined by
\beq \bx^\eta=\nabla_x \eta =\frac{ z_\eta}{\mu  z_\eta \ol z_\eta} e_0 , \ \bx^\theta=\nabla_x \theta =\frac{ z_\theta}{\mu z_\theta \ol z_\theta}e_\theta, \  \bx^\varphi=\nabla_x \varphi =  \frac{1}{\mu z_\varphi e_0}. \label{reciptanvec} \eeq
It is easy to show that $(\bx^\eta)^2=(\bx^\theta)^2=\bx_\eta^{-2}=\bx_\theta^{-2}$, and
\[ z\ol z= \frac{1}{2}(\cos 2 \theta+ \cosh 2 \eta), \ z_\eta \ol z_\eta = \frac{1}{2}(-\cos 2 \theta+ \cosh 2 \eta)=z_\theta \ol z_\theta, \]
and $ z_\varphi \ol z_\varphi =\sinh\eta^2\sin \theta^2 $. We also have
\beq \nabla_x^2 \eta =\frac{\coth \eta}{\mu^2  z_\eta \ol z_\eta} , \ \nabla_x^2 \theta =\frac{\cot \theta}{\mu^2 z_\theta \ol z_\theta}e_\theta, \  \nabla_x^2 \varphi =  0, \label{nabla2recipvec} \eeq
which will be use later. 
 
 For the {\it oblate orthogonal tangent vectors}, $\{\by_\eta,\by_\theta, \by_\varphi \}$, and the corresponding orthogonal reciprocal frame  $\{\by^\eta,\by^\theta, \by^\varphi \}$,
\[ \by_\eta := \partial_\eta \by = \mu z e_0,\  \by_\theta:=  \partial_\theta \by := \mu I_p z e_0  ,\  \by_\varphi:=\mu  z_{\varphi\eta}e_0 = \mu \dot e_p \cosh \eta \sin \theta,   \] 
\[ \by^\eta=\nabla_y \eta =\frac{1}{\mu \ol z} e_0 , \ \by^\theta=\nabla_y \theta =\frac{I_p}{\mu \ol z}e_0, \  \by^\varphi=\frac{1}{\mu \ol z_{\varphi \eta}}e_0 =  \frac{\dot e_p}{\mu \cosh \eta \sin \theta }. \]
  We also have
 \beq \nabla_y^2 \eta =\frac{\tanh \eta}{\mu^2  z \ol z} ,\ \tanh \eta =\frac{z_\varphi}{z_{\varphi \eta}},    \ \nabla_y^2 \theta =\frac{\cot \theta}{\mu^2 z \ol z}, \ \tan \theta =\frac{z_\varphi}{z_{\varphi \eta}}, \  \nabla_y^2 \varphi =  0. \label{nabla2yrecipvec} \eeq
 \section{Spheroidal-graphic projection}   
\begin{figure}[h]
   	\begin{center}  	\includegraphics[width=0.7\linewidth]{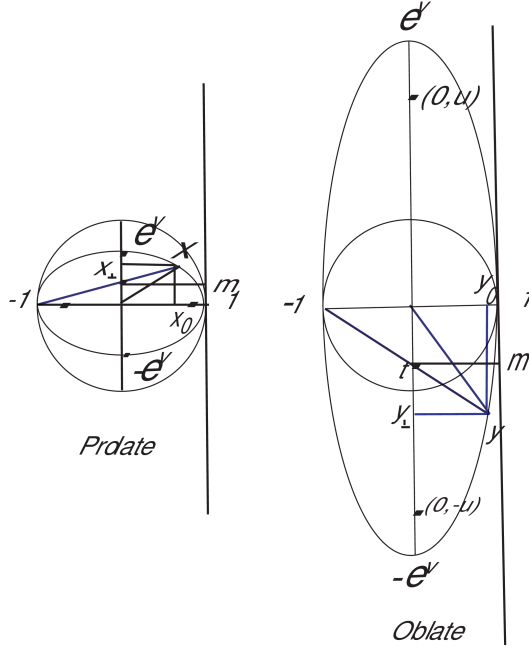}  \end{center}
   	\caption{The elliptical sections of the circumscribed unit prolate Case 3 and inscribed unit oblate Case 1. For $\nu \ge 0$, $e^{2 \nu}-\mu^2 =1$ for Cases 1, 2, and
   		$e^{-\nu}+\mu^2 = 1$ for Cases 3, 4 in Figure \ref{prolateoblate}, respectively. When $\mu \to 0$ and $\nu \to 0$, the prolate and oblate spheroids go to the $3$-sphere.}
   	\label{prolatoblate3}
   \end{figure}

   We now define spheroidal-graphic projection from the point
   $-\be_0$ on the bounding prolate and oblate unit spheroids 1 and 3 in Figure \ref{prolateoblate}, respectively, to the corresponding vertical $(0,x_1,x_2)$, $(0,y_1,y_2)$ planes, shown in Figure \ref{prolatoblate3} as vertical lines. Clearly, as the point $\bx$ moves along the surface of the unit prolate, the projected point $t\be_p:= s \bx_p$  moves in the interior of the disk bounded by the circle with the points $ -e^{-\nu} \be_p$ and $e^{-\nu}\be_p $ on its diameter. Similarly, as the point $\by$ moves along the surface of the unit oblate spheroid, the projected point $t\be_p:= s\by_p$ moves in the interior of the disk bounded by the circle with the points $ -e^{\nu} \be_p$ and $e^{\nu} \be_p$ on its diameter.      
   
   The spheroidal-graphic projections $t \be_p$ for unit prolate and oblate spheroids are easily defined. We have $t=sx_p$ and $t=sy_p$ for
   \beq t = e^{-\nu}\sqrt{\frac{1-x_0}{1+x_0}} \quad {\rm and} \quad t = e^{\nu}\sqrt{\frac{1-y_0}{1+y_0}}\label{mappingprob} \eeq
   in the prolate and oblate cases, respectively. Letting 
   $ \bm = t\be_p  + \be_0 = s(\bx +\be_0)$,
   \beq s=\frac{|t \be_p  + \be_0|}{|\bx + \be_0|}=\frac{1}{x_0+1} \ \iff \ \frac{t \be_p  + \be_0}{1} =\frac{\bx +\be_0}{x_0+1}\ \iff \ x_0 = \frac{\bx - t \be_p  }{t\be_p   + \be_0 }   \label{mappingextra} \eeq 
   in the prolate case, the mapping (\ref{mappingprob}) relating similar triangles reduces to
   \beq    t\be_p   = \frac{\bx +\be_0}{x_0+1}-\be_0  \ \ \iff \ \ t \be_p  =\frac{\bx - x_0 \be_0 }{x_0+1}=\frac{x_1 \be_1+x_2 \be_2}{x_0+1},   \label{mapping2} \eeq
   which implies that $t =\frac{x_1 \cos \varphi + x_2 \sin \varphi}{x_0+1}$. Exchanging $x$'s for $y$'s give the similar result in the oblate case.

   Two important relationships for the both the oblate/prolate case are
   \beq   \mu^2 = \frac{1-\bx^2}{1-x_0^2} \ \ \iff \ \ e^{\pm 2\nu} = 1\pm \mu^2 = \frac{\bx^2-x_0^2}{1-x_0^2}\ \ {\rm and} \ \ x_0 = \frac{e^{\pm 2\nu}-t^2}{e^{\pm 2\nu}+t^2} ,            \label{twoimport} \eeq
   where the $``+"$ sign is chosen for the $\by$-oblate case 1, and the $``-"$ sign is chosen in the $\bx$-prolate case 3, shown Figure \ref{prolatoblate3}. 
   Using the last relationship, we can easily invert the mapping in  (\ref{mappingextra}) or (\ref{mapping2}), in both the oblate-prolate cases,
   getting
   \beq \bx = \frac{2e^{\pm 2\nu}t \be_p +(e^{\pm 2\nu}-t^2)\be_0 }{e^{\pm 2\nu}+t^2} \ \ \iff \ \ \bx+\be_0 = \frac{2e^{\pm2\nu}(t \be_p +\be_0)} {e^{\pm 2\nu}+t^2}.   \label{mapping3}\eeq
   There is an interesting relationship between spheroidal-graphic projection and the Vekua system of equilibrium equations in a spherical shell \cite{Zhg1981}, which will be explored elsewhere.
   
   In both the oblate-prolate cases, when $\mu \to 0$,  $\nu \to 0$, $t=\sqrt{\frac{1-x_0}{1+x_0}}$, the mappings (\ref{mappingprob}) and
   (\ref{mapping3}) go to stereographic
   projection $t \be_p$ from the point $-\be_0$ to a point in the plane of the bivector $\be_{12}$ passing through the origin, 
   \beq t \be_p  + \be_0 =\frac{\bx + \be_0 }{x_0+1} =\frac{(\bx + \be_0)^2 }{(x_0+1)(\bx+\be_0 )} = \frac{2}{\bx +\be_0}, \label{stereop} \eeq
   with the stereographic inverse mapping
   \[ \bx = \frac{2 t \be_p  +(1-t^2)\be_0 }{1+t^2} \ \ \iff \ \ \bx+\be_0 = \frac{2}{t \be_p +\be_0 } . \]
   Stereographic projection has been extensively studied in geometric algebra in
   \cite[pp.111-120]{Sob2015} and \cite{Sob2019}.
   
   The relationships (\ref{mappingprob}) and (\ref{stereop}) can easily be express in spheroidal coordinates in both the oblate-prolate cases 1 and 3 in Figure \ref{prolateoblate}. Since we are assuming that for a fixed $\mu$, $\cosh \eta = \frac{1}{\mu}$ in equation (\ref{mappingprob}), the spheroidal coordinate form of equation (\ref{mapping3}) in terms of $(\eta, \theta, \varphi)$ is 
   \[ t \be_p  =  e^{\pm\nu}\sqrt{\frac{1-x_0}{1+x_0}}\, \be_p =
   \tanh \eta \sqrt{\frac{1- \cos\theta }{1+ \cos \theta}}\,  \be_p   \]
   for $ \be_p = \be_1\cos \varphi  +  \be_2\sin\varphi$.              
   
    \section{Spheroidal gradient and Laplacian}
    
    In the terms of rectangular coordinates 
    \[  \bx = x_0 e_0 +x_1 e_1+x_2 e_2, \quad \by =y_0 e_0 +y_1 e_1+y_2 e_2,  \]
  the gradient and Laplacian take the usual forms
  \[  \nabla_x = e_0 \partial_{x_0} + e_1 \partial_{x_1}+ e_2 \partial_{x_2}, \quad   \nabla_x^2 =  \partial_{x_0}^2 + \partial_{x_1}^2+  \partial_{x_2}^2\]  
  and
   \[  \nabla_y = e_0 \partial_{y_0} + e_1 \partial_{y_1}+ e_2 \partial_{y_2}, \quad   \nabla_y^2 =  \partial_{y_0}^2 + \partial_{y_1}^2+  \partial_{y_2}^2, \]
   respectively.
   
   In prolate spheroidal coordinates, the gradient and Laplacian are respectively given by
   \[ \nabla_x = \frac{e_0}{\mu}z_\eta^{-1} \bigg( \partial_\eta
       -I_p\partial_\theta +z_\eta z_\varphi^{-1}\partial_\varphi \bigg) \]
        \[=  \frac{e_0}{\mu}z_\eta^{-1} \bigg( \partial_\eta
        -I_p\partial_\theta - \Big(J_p\cot \theta +K_p\coth \eta  \Big)\partial_\varphi \bigg)  ,\]
       \[  \nabla_x^2=\frac{1}{\mu^2} \Big(\frac{1}{\ol z_\eta}  \partial_{\eta} + \frac{1}{\ol z_\theta} \partial_{\theta}+\frac{1}{\ol z_\varphi}   \partial_{\varphi}\Big) \Big(\frac{1}{z_\eta}  \partial_{\eta} + \frac{1}{ z_\theta} \partial_{\theta}+\frac{1}{z_\varphi}\partial_\varphi \Big) \]
       \beq = (\nabla_x \eta)^2 \bigg( \partial_\eta^2 +\partial_\theta^2 +\frac{(\nabla_x \varphi)^2}{(\nabla_x \eta)^2}\partial_\varphi^2+\frac{(\nabla_x^2\eta)\partial_\eta +(\nabla_x^2\theta)\partial_\theta }{(\nabla_x \eta)^2}\bigg),   \label{gandlpro} \eeq
   \cite{nested}.    
  In terms of the quaternion $z$, the Laplacian takes the form
  \[ \nabla_x^2=\frac{1}{\mu^2 z_\eta \ol z_\eta}\bigg(\partial_{\eta}^2 + \partial_{\theta}^2+\frac{z_\eta \ol z_\eta}{z_\varphi \ol z_\varphi}\partial_{\varphi}^2+ \coth \eta \, \partial_{\eta} + \cot \theta \,\partial_{\theta} \bigg)    \]     
   \beq =\frac{1}{\mu^2 z_\eta \ol z_\eta}\bigg(\partial_{\eta}^2 + \partial_{\theta}^2+\Big(\cot^2 \theta + \coth^2 \eta\Big)\partial_{\varphi}^2+ \coth \eta \, \partial_{\eta} + \cot \theta \,\partial_{\theta} \bigg),    \label{prolate-laplace} \eeq
   equivalent to the same equation found in \cite[p.\,411]{Hob1931}.
       
   In oblate spheroidal coordinates, the gradient and Laplacian are respectively given by    
         \[  \nabla_y = \frac{e_0}{\mu}z^{-1} \bigg( \partial_\eta
       -I_p\partial_\theta +z z_{\varphi \eta}^{-1} \partial_{\varphi}\bigg) \]
       \[ =  \frac{e_0}{\mu}z^{-1} \bigg( \partial_\eta
       -I_p\partial_\theta - \Big(J_p\cot \theta +K_p\tanh \eta  \Big)\partial_\varphi \bigg)    \]
     \[  \nabla_y^2 =\frac{1}{\mu^2} \bigg(\frac{1}{z_\eta}  \partial_{\eta} + \frac{1}{ z_\theta} \partial_{\theta}+\frac{1}{z_{\varphi\eta}} \partial_{\varphi}\bigg) \bigg(\frac{1}{\ol z_\eta}  \partial_{\eta} + \frac{1}{\ol z_\theta}   \partial_{\theta} +\frac{1}{ \ol z_{\varphi\eta}} \partial_{\varphi}\bigg) \]
     \beq =  (\nabla_y \eta)^2 \bigg( \partial_\eta^2 +\partial_\theta^2 +\frac{(\nabla_y \varphi)^2}{(\nabla_y \eta)^2}\partial_\varphi^2+\frac{(\nabla_y^2\eta)\partial_\eta +(\nabla_y^2\theta)\partial_\theta }{(\nabla_y \eta)^2}\bigg)    .  \label{gandlob} \eeq
     In terms of the quaternion $z$, the Laplacian takes the form
    \[ \nabla_y^2=\frac{1}{\mu^2 z \ol z}\bigg(\partial_{\eta}^2 + \partial_{\theta}^2+\frac{z \ol z}{z_{\varphi\eta} \ol z_{\varphi\eta}}\partial_{\varphi}^2+ \tanh \eta \, \partial_{\eta} + \cot \theta \,\partial_{\theta} \bigg)    \]     
      \beq =\frac{1}{\mu^2 z \ol z}\bigg(\partial_{\eta}^2 + \partial_{\theta}^2+\Big(\cot^2 \theta + \tanh^2 \eta\Big)\partial_{\varphi}^2+ \tanh \eta \, \partial_{\eta} + \cot \theta \,\partial_{\theta} \bigg).  \label{prolate-laplace2} \eeq   
           Note that $(\nabla_x \eta)^2=(\nabla_x \theta)^2$ and $(\nabla_y \eta)^2=(\nabla_y \theta)^2$ in the expressions  (\ref{gandlpro}) and (\ref{gandlob}) above, and that the expressions are the same except for the gradients employed with respect to $\bx $ and $\by $, respectively, \cite{nested}. 
           
       \section{Quaternion gradient and Laplacian} 
     
     Both the prolate and oblate gradients and Laplacians can be expressed in terms of a more fundamental quaternion gradient and Laplacian, as is explored in this section.
     
     Beginning with the results given in (\ref{prolate-laplace}) and (\ref{prolate-laplace2}), the {\it quaternion gradient} is defined by
     \beq \nabla_z:=  \bigg( \frac{1}{z_\eta} \partial_{\eta}+ \frac{1}{z_\theta} \partial_{\theta} +  \frac{1}{z_\varphi} \partial_{\varphi}\bigg)=\frac{1}{z_\eta } \bigg( \partial_{\eta}-I_p\partial_{\theta} +  z_\eta z_\varphi^{-1} \partial_{\varphi}\bigg), \label{quaterniong} \eeq
     and 
     \beq \nabla_{z_\eta}:= \frac{1}{z } \bigg( \partial_{\eta}-I_p\partial_{\theta} +  z z_{\varphi \eta}^{-1} \partial_{\varphi}\bigg). \label{quaterniongeta} \eeq
     Note in the above definitions 
     \[  \frac{1}{z_\varphi}=  \frac{1}{\dot I_p\sinh \eta \sin \theta} =  -\dot I_p  \frac{1}{\sinh \eta \sin \theta}, \]
     $\nabla_x = \frac{e_0}{\mu}\nabla_z$, $\nabla_x^2 =  \frac{1}{\mu^2}\nabla_z \nabla_{\ol z}$, and $\nabla_y = \frac{e_0}{\mu}\nabla_{z_\eta} $, $\nabla_y^2 =  \frac{1}{\mu^2}\nabla_{z_\eta}  \nabla_{\ol z_{\eta}}$,
     where      
     \[ \nabla_{\ol z}:=e_0 \nabla_z e_0 = \bigg( \frac{1}{\ol z_\eta} \partial_{\eta}+ \frac{1}{\ol z_\theta} \partial_{\theta} +  \frac{1}{\ol z_\varphi} \partial_{\varphi}\bigg), \ \ \nabla_{\ol z_\eta}:=e_0 \nabla_{z_\eta} e_0  . \]
     The {\it prolate quaternion Laplacian} is given by
     \beq \nabla_{\ol z}\nabla_z =e_0 \nabla_z e_0 \nabla_z\equiv \mu^2\nabla_x^2=\nabla_z \nabla_{\ol z},                          \label{quaternionL}  \eeq
     and similarly for the  {\it oblate quaternion Laplacian}.
     The quaternion Laplacians are, up to a scalar factor, equivalent to the prolate and oblate Laplacians $\nabla_x^2$ and $\nabla_y^2$ given in (\ref{gandlpro}) and  (\ref{gandlob}), respectively. 
     
     Below is a Table of useful identities:
     
     \begin{itemize}
     	\item[1.] $\nabla_z z = 3, \ \nabla_z \ol z=-1 $.
     	\item[2.] $z \ol z =(\cosh 2 \eta + \cos 2\theta), z_\eta \ol z_\eta =\frac{1}{2}(\cosh 2 \eta - \cos 2\theta), z_\varphi \ol z_\varphi = \sinh^2 \eta \sin^2 \eta$.  
     	\item[3.] $\nabla_x z \ol z= \frac{2 \bx}{\mu^2}, \ \nabla_z z_\eta \ol z_\eta=0 = \nabla_{\ol z} z_\eta \ol z_\eta $. 
     	\item[4.] $z_\eta \ol z + z \ol z_\eta = \sinh 2\eta = -I_p(z_\theta \ol z - z \ol z_\theta ) $,
     	
     	$z_\theta \ol z - z \ol z_\theta =- \sin 2 \theta=I_p(z_\eta \ol z - z \ol z_\eta )  $.
     	\item[5.] $\nabla_x^2 \eta =\frac{\coth \eta}{\mu^2 z_\eta \ol z_\eta}, \ \nabla_x^2 \theta =\frac{\cot \theta}{\mu^2 z_\theta \ol z_\theta}, \ (\nabla_x \eta)^2 =\frac{1}{\mu^2 z_\eta \ol z_\eta}= \frac{1}{\mu^2 z_\theta \ol z_\theta}=(\nabla_x \theta)^2$.
     	
     	\item[6.] $\nabla_y^2 \eta =\frac{\tanh \eta}{\mu^2 z \ol z}, \ \nabla_y^2 \theta =\frac{\cot \theta}{\mu^2 z \ol z}, \ (\nabla_y \eta)^2 =\frac{1}{\mu^2 z \ol z}=(\nabla_y \theta)^2, \ \frac{\nabla_y^2 \eta}{(\nabla_y \eta)^2}=\tanh \eta$.
     	\item[7.] $\frac{z-\ol z}{z+\ol z}=I_p \tanh \eta \tan \theta , \ \frac{z_\eta-\ol z_\eta}{z_\eta+\ol z_\eta}=I_p \coth \eta \tan \theta, \ \frac{|\nabla_y \varphi|^2}{|\nabla_x \varphi|^2 }=\tanh^2 \eta$.
     \end{itemize}
     
     The properties of the quaternions $I_p:=e_p e_0, J_p :=\dot e_p e_0,$ and $K_p := I_p J_p = \dot e_p e_p$, are given below:
     
     \begin{itemize}
     	\item[1.] $I_p^2=J_p^2=K_p^2=-1, \ I_p J_p K_p = -1$.
     	\item[2.] $\dot I_p:=\partial_\varphi I_p = J_p, \ \dot J_p=\partial_\varphi^2 I_p=-I_p, \ \dot K_p =0. $ 
     \end{itemize}
     The fact that $\dot K_p = 0$ is a consequence of $\partial_\varphi \dot e_p = -e_p$.
     
     Clearly the gradients $\nabla_x$, $\nabla_y$, and $\nabla_z$ and $\nabla_{z_\eta}$, are all closely related, since
     \[ \nabla_x =\frac{e_0}{\mu}z_\eta^{-1}  \bigg( \partial_{\eta}-I_p\partial_{\theta} +  z_\eta z_\varphi^{-1} \partial_{\varphi}\bigg)=\frac{e_0}{\mu} \nabla_z     \] 
     and
     \[\nabla_y=\frac{e_0}{\mu} z^{-1} \bigg( \partial_{\eta}-I_p\partial_{\theta} +  z z_{\varphi \eta}^{-1} \partial_{\varphi}\bigg) =\frac{e_0}{\mu} \nabla_{z_\eta}.\]  
     
     \section{Spheroidal solutions to the Laplace equation}
     
      Since prolate an oblate coordinates are one of the 11 systems in which the Laplace equation is separable, harmonic solutions of the equations (\ref{gandlpro}) and (\ref{gandlob}) have the form
    \beq U(\eta,\theta,\varphi)={\cal N}(\eta)\Theta(\theta)\Phi (\varphi) ,  \label{harmolnicsolnp} \eeq
    where $\{{\cal N}(\eta),\Theta(\theta),\Phi(\varphi)\}\in \R$. 
    In the prolate case, separating (\ref{prolate-laplace}) leads to the differential equations,   
    \beq \frac{{d^2\cal N}}{d\eta^2} +\coth \eta  \frac{d\cal N}{d\eta} +\bigg[ -m^2 \coth^2 \eta+n\bigg] {\cal N}=0,  \label{difeqnNp}\eeq
    \beq \frac{{d^2\Theta}}{d\theta^2} +\cot \theta  \frac{d\Theta}{d\theta}- \bigg[ n+m^2 \cot^2 \theta\bigg] {\Theta}=0,  \label{difeqnthetap}\eeq    
    and
    \beq \frac{d^2 \Phi}{d \varphi^2}+m^2 \Phi =0. \label{difthirdp} \eeq
   Separating (\ref{prolate-laplace2}) in the oblate case, only the first equation (\ref{difeqnNp}) changes to 
    \beq \frac{{d^2\cal N}}{d\eta^2} +\tanh \eta  \frac{d\cal N}{d\eta} +\bigg[ -m^2 \tanh^2 \eta+n\bigg] {\cal N}=0,  \label{difeqnNo}\eeq
    the other two equations (\ref{difeqnthetap}) and (\ref{difthirdp}) remaining the same. Solutions involving hypergeometric functions \cite{HB1969} are shown in the Mathematica Package in the Appendix. However, equivalent but much more compact and workable solutions have been found in terms of Legendre functions of the first and second kind, see \cite[p.\,47]{BKM1976} and \cite[pp.\,413,422]{Hob1931}. An extensive discussion of the issues involved in the solutions of the Helmholtz and Laplace equations in terms of their associated   Lie algebras and symmetry groups is given in \cite[pp.\,36-43]{BKM1976}, \cite{Miller1974}.
       
     Following Garabedian \cite{PGar1953}, and Hobson, \cite[p.422]{Hob1931}, the second order differential equations  have the respective interior/exterior harmonic spheroidal solutions of the form
    \beq P_{n,m}[\cos \theta]P_{n,m}[\cosh \eta]\pmatrix{\cos \cr \sin}(m\varphi),      \label{inU} \eeq 
    and 
   \beq  P_{n,m}(\cos \theta)Q_{n,m}(\cosh \eta)\pmatrix{\cos  \cr \sin}(m\varphi),      \label{exU} \eeq 
      respectively, where $P_{n,m}$ and $Q_{n,m}$ are symbols for the respective Legendre Polynomials of the first and second kind \cite{wofram}, and where 
    \[  C_m[\alpha]:=\cos( m\varphi)=m \sum_{k=0}^{[\frac{m}{2}]}(-1)^k\frac{(m-k-1)!}{k!(m-2k)!} 2^{m-2k-1}\alpha^{m-2k} . \] 
        
    From the prolate and oblate cases (\ref{pomega}) - (\ref{obomega}) involving $\bx$ and and $\by$, by substituting expressions for $\cosh \eta$ and $\cos \theta$, using Hobson's solutions (\ref{inU}) and (\ref{exU}), we get harmonic polynomial solutions in terms of the variables $\{x_0,x_1,x_2\}$ in the prolate case, and $\{ y_0,y_1,y_2\}$ in the oblate case,   \cite[pps.\,413,\,422]{Hob1931}. The theoretical framework for the study of different separable solutions is considered in the next Section. 
    
   \section{Lie algebra ${\cal E}(3)$ of symmetry operators}
   As explained in \cite[pp.36-43]{BKM1976}, the six dimensional real {\it Lie algebra} ${\cal E}(3)$ of the {\it Euclidean symmetry group} $E(3)$ is generated by a basis of six symmetry operators
   \beq {\cal J}_k = \be_k \cdot {\cal J}, \ {\cal P}_k = \partial_k, \label{oldsymop} \eeq
   for ${\cal J}_x:=-\bx \times \nabla_x$ and $k=0,1,2$.   
   The basic theory of this Lie algebra is developed here in a new way utilizing the rich structure of the geometric algebra $\G_3$.
   
   Let $\ba , \bb \in \R^3$ be arbitrary constant vectors in $\G_3^1$. Define
   the {\it scalar operator} $P_a$ and the {\it vector operator} $J_b$ by
    \beq P_a := \ba \cdot \nabla_x, \quad J_b := \bb \w \bx \w \nabla_x,
         \label{newsymop}  \eeq  
    and  
    \[ \bJ_x:= i\, \bx \w \nabla_x=- \bx \times \nabla_x = {\cal J}_x,   \]
   for $i := \be_{012}$. The interesting relationship 
   \beq \bJ_x^2 = (i\, \bx \w \nabla_x)^2= \bx^2 - \bx \cdot \nabla_x -(\bx \cdot \nabla)^2,     \label{JxJx} \eeq
   follows after a rather tricky calculation.  
      
   The close relationship between the definitions (\ref{oldsymop}) and (\ref{newsymop}) is easily found,
   \[ {\cal P}_k= P_{e_k}=\be_k \cdot \nabla_x, \ \ {\rm and} \ \ {\cal J}_k =\be_k \cdot \bJ_x = i\, \be_k \w \bx \w \nabla_x  .     \]
   We can now state the basic Lie algebra bracket relationships among the symmetry operators:
  \beq [P_a,P_b]= 0, \ [ J_a,P_a]=-i P_{a\times b}, \ [J_a,J_b]= i J_{a\times b}. 
   \label{liealgrel} \eeq
   By the symmetry Lie algebra $\cal S$ of symmetry operators $S_{\ba,\bb}=P_{\ba}+J_{\bb}$, we mean
   \beq {\cal S}:=\{S_{\ba,\bb}|\ \ \ba, \bb \in \R^3 \} \label{symalgebraS} \eeq
   Thus, a general symmetry operator $S_{\ba,\bb}$ is the sum of a scalar and pseudo- scalar operator parts. Since $i=\be_{012}$ is in the center $\cal Z$ of $\G_3$, a symmetry operator will naturally commute with any constant multivector in $\G_3$. 
   
   The importance of the symmetry Lie algebra $\cal S$ follows from the fact that the subset of symmetry operators ${\cal L}\subset {\cal S}$, with the property that  $S_{\ba,\bb}\Psi$ is a solution of the Laplace or Helmholtz equation whenever
   $\Psi$ is an analytic solution, make up a Lie sub-algebra of $\cal S$, \cite[p.36]{BKM1976}, \cite{Miller1974}. Furthermore, as noted by these authors, each of these 11 systems of orthogonal coordinates systems in which the Helmholtz equation separates corresponds to a pair of commuting second order operators in the enveloping algebra of ${\cal E}(3)$ of $\cal L$. Studying properties of the Lie algebra $\cal L$, of the Helmholtz equation, for example, gives insight into how the hypergeometric solutions to the prolate and oblate Laplace equations (\ref{gandlpro}) and (\ref{gandlob}) are related to the equivalent famous solutions given by the Legendre polynomial solutions (\ref{inU}) and (\ref{exU}).

   \section{Geometric analysis verses Clifford analysis}
    Clifford analysis \cite{DSS1992} is laid down in terms of the more comprehensive geometric analysis, and in such a way that it is easy to translate any equation in Clifford analysis into its equivalent expression in the geometric analysis, and vice-versa \cite{H/S,SNF,DSS1992}. Applications and examples are given.

    Let $\bx \in \G_n^1$ be the real {\it position vector} in the geometric algebra
    $\G_{n+1}:=\R[\be_0,\be_1, \ldots , \be_n]$ of Euclidean space $\R^n$. Thus,
    \[ \bx = \sum_{k=0}^n x_k \be_k = (x_0, x_1, \ldots, x_n)\in \R^{n+1} .  \] 
    To get the equivalent paravector $\bX \in \G_{0,n}^{0+1}$, write
    \[ \bX := \bx \be_0 = \bx \cdot \be_0 + \bx\w \be_0 =x_0 + \ul \bX = (x_0, x_1, \ldots, x_n)\in \R^{n+1}  ,   \] 
    where 
    \[\bx \w \be_0 := \ul\bx \be_0=    \ul \bX := \sum_{k=1}^n x_i \be_{k0} \in \G_{n+1}^{2}\subset \G_{n+1}^+ \widetilde{=}\G_{0,n}. \]
       Also defined in Clifford analysis is the {\it complex conjugate} 
    $ \ol{\bX} := x_0 - \ul \bX=\be_0 \bx $. Clearly
    \[  \bX = \bx \be_0=\be_0 \ol \bx  \quad \iff \quad \bx = \bX \be_0 =\be_0 \ol X,     \] 
   or equivalently,
      \[ \ol X = \be_0 X \be_0   , \ \ {\rm and} \ \ \ol \bx := \be_0 \bx \be_0=\ol X \be_0  .\] 
     
     In the geometric algebra $\G_{n+1}$ the dot and wedge product are simply defined by
     \[ \ba \bb = \frac{1}{2}(\ba \bb + \bb \ba)+  \frac{1}{2}(\ba \bb - \bb \ba) \equiv \ba \cdot \bb + \ba \w \bb, \]
     the symmetric part being the dot-product and the anti-symmetric part the wedge-product. To see how this carries over to Clifford analysis, write
     \[ \ba \cdot \bb = \frac{1}{2}(\ba \bb + \bb \ba) = \frac{1}{2}(\ba \be_0 \be_0   \bb + \bb \be_0 \be_0   \ba)=\frac{1}{2}\big(\bA \ol \bB + \bB \ol \bA\Big),   \]
and similarly
 \[ \ba \w \bb = \frac{1}{2}(\ba \bb - \bb \ba) = \frac{1}{2}(\ba \be_0 \be_0   \bb - \bb \be_0 \be_0   \ba)=\frac{1}{2}\big(\bA \ol \bB - \bB \ol \bA\Big)   \]
 for the outer product.
    
    To see that the even sub-algebra $\G_{n+1}^+$ is isomorphic to the geometric
    algebra $\G_{0,n}$, following \cite[pp.\,65-80]{Sob2019}, write
    \[ \G_{0,n} :=\R[f_1,\ldots, f_n] \widetilde = \R[e_{10},\ldots,e_{n0}] \widetilde =\G_{n+1}^+ \subset \R[e_0, e_{10},\ldots,e_{n0}] = \G_{n+1}.        \] 
    It is now an easy exercise to translate any expression in Clifford analysis to a corresponding expression in geometric analysis as laid down in \cite{H/S,SNF}. Thus for the {\it Clifford paravector} $\bX = \bx \be_0$,
    \[ \bX \ol \bX = \bx \be_0 \be_0 \bx = \bx^2= \sum_{k=0}^n x_k^2.      \] 
   
    In Clifford analysis, the operator $\partial_X$ and $\partial_{\ol X}$ are defined by
    \[  \partial_X :=\partial_0+ \sum_{k=1}^{n}\be_{k0}\, \partial_k, \quad {\rm and} \quad \partial_{\ol X} :=\partial_0- \sum_{k=1}^{n}\be_{k0}\, \partial_k   \]
which translate to
   \[  \nabla_{\bx}=\partial_{X}e_0 = e_0 \partial_{\ol X} .    \] 
   It follows that the Laplacian $\nabla_{\bx}^2 =\nabla_{\bx}\be_0 \be_0 \nabla_{\bx} \equiv \partial_{X}\partial_{\ol X}$.          
   
  One important application is the so called Cauchy-Kovalevska (CK) extension, which is a construction of a higher order monogenic function from a given monogenic function \cite{CFM2017}, \cite[p.\,151]{DSS1992}. Following \cite[(3.2)]{CFM2017}, as a simple example of the CK extension, consider
  \beq CK[(\ul\bx \be_0  )^k]=CK[-\ul\bx^k  ] :=-\sum_{k=1}^{n}(x_k+x_0 \be_{k0})^k.    \label{ckext} \eeq 
   For $k=2=n$ and  $\ul \bx= \bx_p$,   
  \[ CK[\bX_p^2] = 2x_0^2 + (\bx_p\be_0)^2-2x_0 \bx_p \be_0= F(\bX_p)\be_0\be_0 =f(\bx)\be_0    \]
  is monogenic for $f(\bx )=  2x_0^2-\bx_p^2 -2x_0\bx_p \be_0 $.
  Checking,
  \[ \partial_{\ul\bX}F[\bX_p]=0=\be_0 \nabla_{\bx}  F[\bX_p]=\be_0 \nabla_{\bx} f(\bx)\be_0. \]  
  More generally, it is easy to show that $ \nabla_{\bx} f(\bx)=0$ for
  \[ f(\bx):=\sum_{k=1}^n (x_k +x_0 \be_{k0} )^{k_i}.    \]    
   
\begin{figure}[h]
   	\begin{center}  	\includegraphics[width=0.7\linewidth]{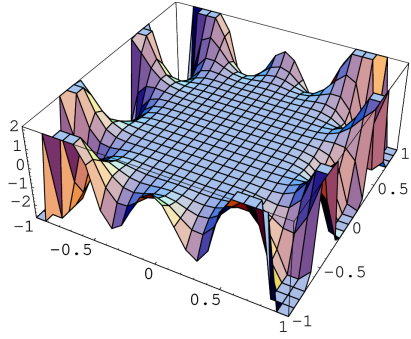}  \end{center}
   	\caption{Shown is the graph of    
   	$\nabla_\bx QM[11]=-11 x_0^{10} +165x_0^8x_p^8 -462x_0^6x_p^4+330 x_0^4x_p^6 -55 x_0^2 x_p^8+ x_p^{10}  $.}
   	\label{rapid}
   \end{figure}

    The idea of a CK extension suggests that the study of {\it quasi-monogenic} (QM) functions $QM[k]$, defined in $\G_3$ by 
   \beq QM[k]:=(x_0 \be_0 - x_p \be_{p0})^k \be_0  , \label{QM-function} \eeq
   is of interest \cite{BCSR1989,Miller1968}. In cylindrical coordinates, the operator $\nabla_\bx$ has the form
   \[ \nabla_{\bx} = (\nabla_{\bx}x_0) \partial_0 +(\nabla_\bx x_p)\partial_p + (\nabla_\bx \varphi)\partial_\varphi = \be_0 \partial_0 + \be_p\partial_p + \frac{\dot\be_p}{x_p}\partial_{\varphi}    ,   \] 
   \cite{nested}. We find that
   \[ \nabla_{\bx} QM[k] = \Big(\be_0 \partial_0  + \be_p\partial_p\Big) \nabla QM[k]+ \frac{\dot\be_p}{x_p}\partial_{\varphi}QM[k] \] \[=\frac{\dot\be_p}{x_p}\partial_{\varphi}QM[k]=\frac{\dot\be_p}{x_p}\cdot \Big(\partial_{\varphi}QM[k]\Big),    \]
   showing that $\nabla_{\bx}\w f(\bx )=0$. Whereas $f(\bx)$ is not monogenic, it is curl-free and $\nabla QM[k]$ rapidly approaches zero in the unit disk in the plane of the bivector $\be_{p0}$, see Figure \ref{rapid}. The CK extension has also been studied in Hermitian Clifford Analysis \cite{BLSS2011}.

   An interesting property of the quasi-monogenic functions (\ref{QM-function}) is that for $k=1,2,3$, the modified functions 
  \[   QM(1)+x_0 \be_0, \ QM(2)+x_0^2 \be_0 ,  \  QM(3)+x_0^3 \be_0 -\frac{x_p^3}{4}\be_p  \]
  are monogenic. Are there other values of $k$ for which
  the function $QM[k]$ can be suitably modified to be monogenic?

   In geometric analysis, the {\it Cauchy kernel} is defined by,
 \beq g(\bx):=\frac{\bx-\by}{|\bx - \by|^{n+1}} =\be_0\frac{\ol\bX-\ol\bY}{|\bX - \bY|^{n+1}} =\frac{\bX-\bY}{|\bX - \bY|^{n+1}}\be_0,    \label{Cauchyk} \eeq
 \cite[p.\,237]{SNF}.
 It is one of the most important examples of a {\it monogenic function}, satisfying 
 \[ \nabla_{\bx} g(\bx) = 0 = \nabla_{\bx}\be_0 \be_0 g(\bx)=\partial_X G(X) , \]       
 where $G(X):= g(\bx) =\frac{\bX-\bY}{|\bX - \bY|^{n+1}}\be_0 $.

    Another interesting method for generating a higher order monogenic functions is by way of the {\it hypercomplex generalized geometric series} of the Cauchy kernel. Starting with the geometric Cauchy kernel function (\ref{Cauchyk}), and employing the complimentary methods of Clifford analysis and geometric analysis, 
  \[  f(\bx)=\frac{\be_0 - \bx}{|\be_0  -\bx|^{n+1}}=( \be_0 -\bx )\Big[(\be_0- \bx )\be_0\be_0(\be_0- \bx  )  \Big]^{-\frac{n+1}{2}}  \]
  \[ =  \be_0 \Big[\be_0 (\be_0- \bx )\Big]^\frac{2}{2}\Big[\be_0(\be_0-\bx)\Big]^{-\frac{n+1}{2}}  \Big[(\be_0-\bx)\be_0\Big]^{-\frac{n+1}{2}}     \]
 \[ = \be_0 \Big[\be_0(\be_0-\bx)\Big]^{-\frac{n-1}{2}}\Big[(\be_0-\bx)\be_0\Big]^{-\frac{n+1}{2}}  \]
 \beq = \be_0 \Big[(1-\ol\bX)\Big]^{-\frac{n-1}{2}}\Big[(1-\bX)\Big]^{-\frac{n+1}{2}}
 =\be_0F(\bX)   ,    \label{hypgeom} \eeq
 where $F(\bX)$ is the Clifford analysis Cauchy kernel \cite[(3.7)]{CFM2017}.
 By expanding each of the expressions defining $F(\bX)$ in a binomial series, the authors' of this reference obtain many beautiful results of hypercomplex generalized geometric series.	
\section*{Appendix: Mathematica Package (Prolate case)}

 \begin{figure}[h] 
 	\centering
 	\includegraphics[width=0.85\linewidth]{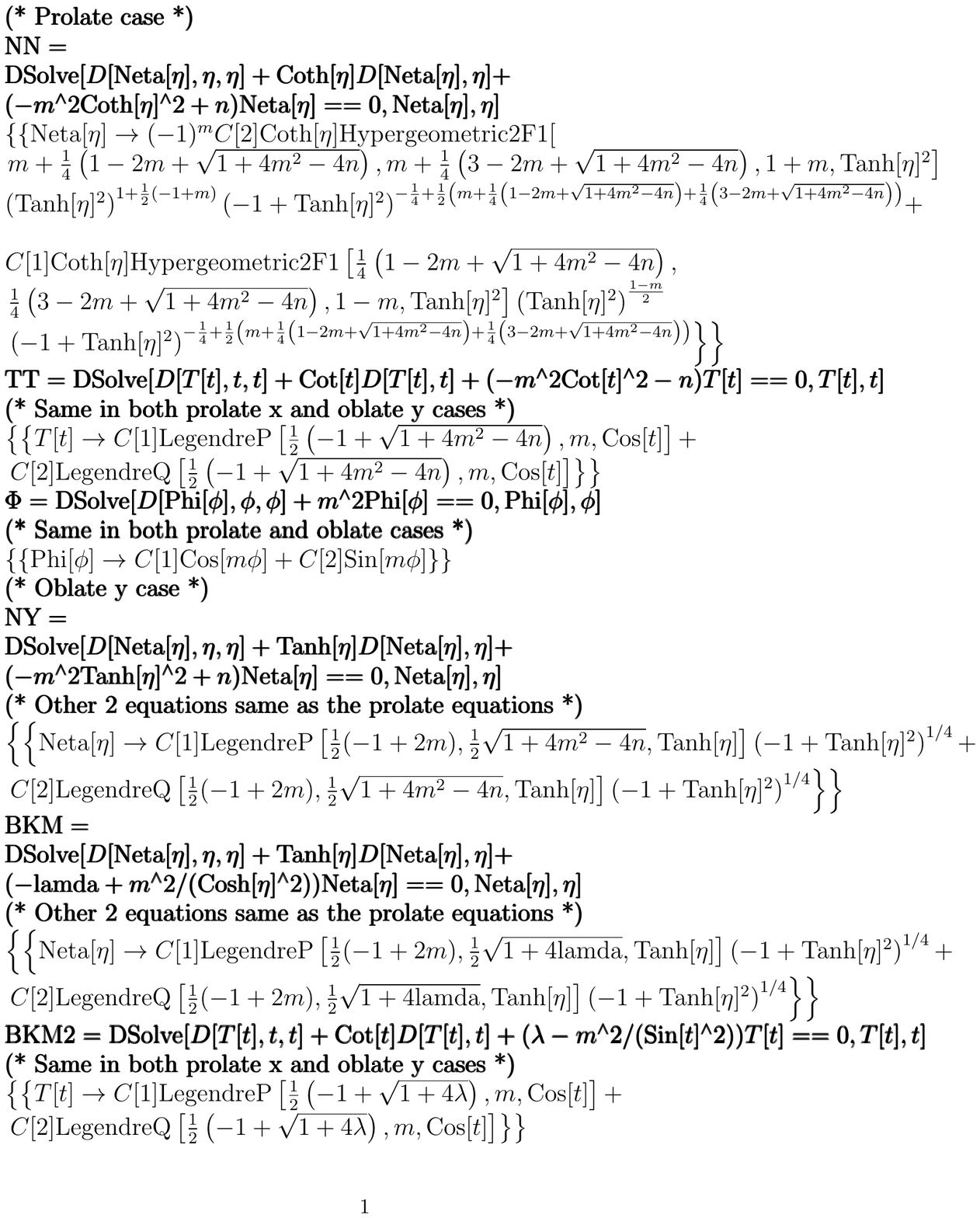}
 	\label{fig:dsolve1}
 \end{figure}
  \section*{Acknowledgement} 
  The author is grateful to Universidad de Las Americas, Puebla for many years of support.

\end{document}